\documentclass[twoside,12pt]{article}
\usepackage{amsmath,amscd,amsthm,amssymb,amsxtra,latexsym,epsfig,epic,graphics}
\def\antiddot{\mathinner{\mkern1mu\raise1pt\vbox{\kern7pt\hbox{.}}\mkern2mu
        \raise4pt\hbox{.}\mkern2mu\raise7pt\hbox{.}\mkern1mu}}

\newcommand{\PP}{{\mathbb P}}

\newcommand{\s}{\mathcal}

\newcommand{\cO}{{\s O}}

\newcommand{\cI}{{\s I}}
\newcommand{\cL}{{\s L}}

\newcommand{\punkt}{\hspace{-.3ex}\raise.15ex\hbox to1ex{\Huge.}}

\def \fix#1 {{\hfill\break \bf (( #1 ))\hfill\break}}

\DeclareMathOperator{\Sym}{Sym}
\DeclareMathOperator{\reg}{reg}

\DeclareMathOperator{\Proj}{Proj}

\DeclareMathOperator{\codim}{codim}

\newcommand{\gm}{\mathfrak m}

\newcommand{\Mac}{{\texttt {MACAULAY2}}}

\newtheorem{theorem}{Theorem}[section]

\newtheorem{proposition}[theorem]{Proposition}
\newtheorem{corollary}[theorem]{Corollary}
\newtheorem{conjecture}{Conjecture}[section]
\theoremstyle{definition}

\date{}
\title{Powers of Ideals and Fibers of Morphisms} 
\author{David Eisenbud and Joe Harris
\footnote{The authors gratefully acknowledge 
partial support by the National Science Foundation
during the preparation of this paper.}}

\begin{document}

\maketitle

\begin{abstract}
Let  $X\subset \PP^n=\PP^n_F$ be a projective scheme
over a field $F$, and let
$\phi:X\to Y$ be a finite morphism. Our main result is a formula
in terms of global data
for the maximum of $\reg \phi^{-1}(y)$,
the Castelnuovo-Mumford regularity
of the fibers of $\phi$ over $y\in Y$ , 
where $\phi^{-1} (y)$ is considered as a subscheme of $\PP^n$.

From an algebraic point of view, our formula is related to the
theorem of 
Cutkosky-Herzog-Trung \cite{CHT1999} and 
Kodiyalam \cite{Kodiyalam2000} 
showing that
for any homogeneous ideal $I$ in a standard graded algebra $S$,
$\reg I^t$ can be written as $ dt+\epsilon$ for some non-negative 
integers $d,\epsilon$ and all large $t$. 
In the special case where $I$ contains a power of $S_+$
and is generated by forms of a single degree, our formula gives an interpretation
of $\epsilon$: it is one less than the maximum 
of $\reg \phi^{-1}(y)$, where $\phi$ is  the morphism associated to $I$.

These formulas have strong consequences for ideals generated by
generic forms.

\end{abstract}

\section*{Introduction}
In this note, all schemes will be
projective over an arbitrary field $F$.
For any projective scheme $X\subset \PP^n$ we write $S_X$ for the
homogeneous coordinate ring of $X$, and $I_X$ for the homogeneous
ideal of $X$. We denote by $\reg X$ the Castelnuovo-Mumford regularity
of $I_X$ (if $X=\PP^n$ we make the convention that $\reg X = 1$). 

If $\phi: X\to Y$
is a finite morphism, then the degree of the fiber
$X_y=\phi^{-1}(y)$ is a semicontinuous function of $y\in Y$, and is thus bounded.
It follows that the Castelnuovo-Mumford regularity of
$X_y$, where $X_y$ considered as a subscheme of $\PP^n$, is also bounded.
Our main result in the form of  Corollary \ref{main1}, gives an algebraic formula for $\max \reg(X_y)$
in terms of global data.

A particularly interesting case occurs when $\phi$ is  a morphism induced by a linear projection
$\phi: \PP^n\to \PP^s$.

\begin{theorem}
\label{geometric sample}
Let $X\subset \PP^n$ be a projective scheme with 
homogeneous coordinate ring $S_X$,
and let $\phi: X\to \PP^s$ be a linear projection
whose center does not meet $X$,
defined by an $s+1$-dimensional vector space of
linear forms $V$. Let
$I\subset S_X$ be the ideal generated by $V$, and let 
$\gm$ be the maximal homogeneous ideal of $S_X$.
The maximum of the Castelnuovo-Mumford regularities
of the fibers of $\phi$ over closed points of $\PP^s$ is one more than
the least $\epsilon$ such that, for large $t$,
$$
\gm^{t+\epsilon} \subset I^t.
$$
\end{theorem}

In the situation of Theorem \ref{geometric sample} the number
$t+\epsilon$ is equal, for large $t$, to the Castelnuovo-Mumford regularity of $I^t$
or the corresponding ideal sheaf (see \S 1). Thus
Theorem \ref{geometric sample} clarifies
the following beautiful result of
Cutkosky-Herzog-Trung \cite{CHT1999}, 
Kodiyalam\cite{Kodiyalam2000}, and 
Trung-Wang \cite{TW2005}, at least in a special case.

\begin{theorem}\label{C et al}
If $I$ is a homogeneous ideal in the polynomial ring 
$S=F[x_0,\dots, x_n]$, and $M$ is a finitely generated graded
module over $S$,  then there are non-negative integers $d,\epsilon$ such that
\[
\reg(I^tM) = dt+\epsilon\qquad \hbox{for all $t\gg 0$}. \eqno{\qed}
\]
\end{theorem}

If
$I$ is generated by forms of a single degree $\delta$ and contains a 
nonzerodivisor on $M$, then $d=\delta$. More generally,
Kodiyalam \cite{Kodiyalam2000} proves that
$d$ is the smallest number $\delta$ such that $I^tM=I_{\leq \delta}I^{t-1}M$ 
for large $t$,
where
 $I_{\leq \delta}$ denotes the ideal generated by the elements of
$I$ having degree at most $\delta$.

By contrast, the value of
$\epsilon$ has been mysterious.
Theorem \ref{geometric sample}
gives an interpretation of $\epsilon$ in a special case.
This seems to be new even for ideals generated in a single degree in
a polynomial ring in 2 variables, where Theorem \ref{geometric sample}
yields
the following.

\begin{corollary}\label{two vars}
Suppose that $I\subset F[x,y]$ is an ideal generated by a vector space
$V$ of forms of degree $d$, and that $F$ is algebraically closed. Assume that
the greatest common divisor of the forms in $V$ is 1.
For $V'\subset V$, let
 $r_{V'}$ be the degree of the greatest common divisor of the forms in $V'$,
 and let 
$$
r:=  \max\{r_{V'} \mid V'\subset V \text{a subspace of codimension 1}\}.
$$
If $t\gg 0$, then
$
\reg I^t = dt+ r-1.
$\qed
\end{corollary}

The corresponding result holds in the case of polynomial rings in more variables,
(and also follows from Theorem \ref{geometric sample}) 
if we assume that $I$ is primary to the maximal homogeneous
ideal and redefine $r_{V'}$ to be the maximal degree in which the local cohomology module
$H^1_\gm(R/(V'R))$ is nonzero (Proposition \ref{easy bound-generalized}). 

In the case of 2 variables we may think of $V$ as defining a morphism
$\PP^1 \to \PP(V)$. When this morphism is birational
 the number $r$ can also be interpreted as the maximum multiplicity of a point on the 
 image curve.

The first author and Roya Beheshti \cite{BE} have conjectured 
 that the regularity
of every fiber of a general linear projection of a smooth 
projective variety $X$ to $\PP^{\dim X +c}$, for $c\geq 1$,
is bounded by $1+(\dim X)/c$. 
Translating this conjecture by means of Theorem \ref{geometric sample},
we get:

\begin{conjecture}[Beheshti-Eisenbud \cite{BE}] 
\label{conjecture 1}
Let $R$ be a standard graded $F$-algebra of dimension $n+1$, and let
$\gm$ be the maximal homogeneous ideal of $R$. Suppose that $R$ 
is a domain with isolated singularity. If $F$ is infinite and
$I\subset R$ is an ideal generated by $n+1+c$ general linear forms,
then 
$$
\gm^{t+\epsilon}\subset I^t \quad \text{\rm for all $t\gg 0$}
$$
holds with $\epsilon =\lfloor n/c \rfloor$.
\end{conjecture}

It is easy to see that
Conjecture \ref{conjecture 1} holds if $c>n$, and it is known to hold in many 
other cases as well (see Beheshti-Eisenbud \cite{BE} for a survey),
This also gives some new 
information about ideals generated by generic forms of higher degree. The following
is a typical example. 
Amazingly, we can say no more  than this even
if we assume that $R$ is the polynomial ring $F[x_0,\dots, x_n]$.

\begin{corollary}
Suppose that $R$ is a standard graded algebra of dimension $n+1$ over a field
of characteristic 0,
and that $R$ has at most an isolated singularity. 
If $I=(f_1,\dots, f_{n+2})$ is an ideal generated $n+2$
generic forms of degree $d$, and $n\leq 14$, then
$$
\gm^{t+n}\subset I^t \quad \text{\rm for all $t\gg 0$} 
$$
\end{corollary}

\begin{proof} The linear series given by $f_1,\dots, f_{n+2}$
defines a generic linear projection of Proj $R$.
 By Mather \cite{mat2}, generic projections in this range of dimensions are
 stable maps. Mather \cite{mat1} gives a local classification of the multigerms
 of such stable maps, from which it follows that the degree of the fibers,
 and thus their regularities, are bounded by $n+1$. The desired formula
 now follows from Theorem \ref{geometric sample}. 
\end{proof}

We do not currently know any function $\epsilon$ of $\dim R$ and $c$ alone that
makes the formula in Conjecture \ref{conjecture 1} true. But there is
an elementary estimate, whose proof we will give in \S 1:

\begin{proposition}
\label{easy bound}
Let $R$ be a standard graded $F$-algebra of dimension $n+1$, and let
$\gm$ be the maximal homogeneous ideal of $R$. If
$I\subset R$ is an ideal generated by linear forms, and if
$R/I$ has finite length,
then 
$$
\gm^{t+\epsilon}\subset I^t \quad \text{\rm for all $t\gg 0$}
$$
holds with $\epsilon = \reg R -1$. 
If $X=\Proj R$ is geometrically reduced and connected in
codimension 1, then the same formula holds with
$\epsilon = \deg X -{\rm codim }\ X$.
\end{proposition}

It was conjectured by the first author and Shiro Goto in \cite{Eisenbud-Goto} that
if $X=\Proj R$ is geometrically reduced and connected in
codimension 1, then
$\reg X \leq \deg X - {\rm codim }\ X +1$,
which would say that the first bound given is always sharper than the second,
as well as more general.

In Section \ref{elementary section} we prove a sharp form
of  Theorem \ref{C et al} in the special case of interest for this paper. 
We also give the
proof of a generalization of Proposition \ref{easy bound}.
Section \ref{products} contains our main result, from 
which we derive Theorem \ref{geometric sample}.

We are grateful to Craig Huneke, with whom we first discussed the
problem of identifying the number $\epsilon$ in 
Theorem \ref{C et al}. After some experiments  using 
Macaulay2 \cite{M2}, he suggested the result in 
Corollary \ref{two vars}, which led us to the results of this paper.

\section{The Regularity of Powers}
\label{elementary section}

Throughout this paper we write $S=F[x_0,\dots, x_n]$ and set $\gm = (x_0,\dots, x_n)$,
the homogeneous maximal ideal.

In the case of most interest for this paper, Theorem \ref{C et al}
can be strengthened as follows. The result also sharpens
Theorem 4 of Chandler \cite{Chandler}, where the line of research
leading to Theorem \ref{C et al} began.

\begin{proposition} 
\label{decreasing}
Let $M$ be a finitely generated graded $S$-module generated in degree 0,
and let $I\subset S$ be a homogeneous ideal generated by forms of degree $d$.
If $M/IM$ has
finite length, but $M$ does not, then we may write
\begin{enumerate}
\item $\reg M/I^tM= dt+f_t-1$, with $f_1\geq f_2\geq \cdots\geq 0$.
\item $\reg I^tM= dt+e_t$, with $e_1\geq e_2\geq \cdots\geq 0$.
\end{enumerate}
Moreover, $e:=\inf\{e_t\}=\inf\{f_t\}$, and we have $\reg I^tM=dt+e$ for $t\gg 0$.
\end{proposition}

\begin{proof}
We first prove the inequalities of part 1. Since $M/I^tM$ has finite length 
the assertion $\reg M/I^tM = dt+f_t-1$ means that $f_t$ is the smallest number such that
 $I^tM$ contains
all the graded components of $M$ with degree $\geq dt+f_t+1$. By our hypotheses
on the degrees of generators of $M$ and $I$, this
is equivalent to the assertion $\gm^{f_t+1}I^tM= \gm^{dt+f_t+1}M$.
A priori we have
$\gm^{f_t+1}I^tM\subset \gm^{d+f_t+1}I^{t-1}M\subset \gm^{dt+f_t+1}M$,
so if $\reg M/I^tM \leq dt+f_t$ then these three terms are all equal.

From these equalities we deduce 
$$
\gm^{f_t+1}I^{t+1}M= I \gm^{f_t+1}I^tM = I \gm^{d+f_t+1}I^{t-1}M = \gm^{d+f_t+1}I^{t}M = \gm^{d(t+1)+f_t+1}M
$$
so $f_{t+1}\leq f_t$. Considering the degrees of the generators of $I$ and $M$ we
see that $f_t+1\geq 0$ for each $t$, completing the proof of part 1.

Turning to the assertion of part 2, it is obvious from the consideration 
of degrees that $e_t\geq 0$. To prove that $e_t\geq e_{t+1}$, 
let $N$ be the largest submodule
of finite length in $M$. If $N=0$, then since $M/I^tM$ has 
finite length, we see from the local cohomology characterization
of regularity that 
$$
\reg I^tM = \max(\reg M, 1+\reg M/I^t M)
$$
so part 2 follows from part 1 in this case. Moreover, since
$\reg M/I^tM$ increases without bounds, we see
that for large $t$ we will have $e_t=f_t$.

We can reduce the general case to the case $N=0$ by considering the
exact sequence
$$
0\to I^tM\cap N \to I^tM \to I^t(M/N)\to 0.
$$
Since $I^tM\cap N$ has finite length, 
while $I^t(M/N)$ has no finite length submodule  except 0,
$$
\reg I^tM = \max(\reg (I^t\cap N), \reg (I^t(M/N))).
$$
If we replace $t$ by $t+1$
the term $\reg (I^t\cap N)$  does not increase, while $\reg (I^t(M/N))$
increases by at most $d$, proving that $e_t\geq e_{t+1}$.
Because $\reg (I^t(M/N))$ grows without
bound, it eventually dominates, and we again get $e_t=f_t$ for large $t$. \end{proof}

We remark that Proposition \ref{decreasing} does not hold if we drop the assumption
that $M/IM$ has finite length. As shown by
Sturmfels \cite{Sturmfels}
it is not true in general that $\reg I^2M\leq \reg I M+d$.
 For example, with $M=S$ and char $F \neq 2$,
the ideal associated to a triangulation of the 
projective plane
 has a linear resolution ($\reg I = 3$), but its square does not ($\reg I^2 = 7 > 2 \times 3$).
 Conca \cite{Conca}, gives examples with $\reg I^n = n\reg I$
 but $\reg I^{n+1} > (n+1)\reg I$ for
arbitrary $n$.

We now turn to the proof of Proposition \ref{easy bound}. The first estimate
is a Corollary of the following result:

\begin{proposition}
\label{easy bound-generalized}
Let $M$ be a graded module of dimension $n$ over a polynomial ring
$S=F[x_0,\dots,x_r]$, and let $I$ be an ideal generated by forms of 
degree $d$ such that $M/IM$ has finite length. For every $t>0$ we have
$$
\reg M/I^tM \leq td+ \reg M+(n-1)(d-1)-1 \qquad \hbox{\rm for every $t>0$.}
$$
Moreover, equality holds when the generators of $I$ form a regular
sequence on $M$.  
\end{proposition}

\begin{proof}
If $I$ is generated by a regular sequence on $M$, then
one can obtain a resolution of $M/I^tM$ by tensoring a resolution
of $M$ with one for $S/I^t$ (obtained, for example, as an Eagon-Northcott
complex) and from this one computes the regularity at once.
(This much does not use the hypothesis that $M/IM$ has finite length.)

When $M/IM$ has finite length, we may begin
by replacing $I$ by a smaller ideal, generated by a system
of parameters of degree $d$ on $M$---in this case,
the
regularity of $M/I^tM$ is simply the degree of
the socle, so it can only increase. It is not hard to give
an elementary argument using induction on $t$. Alternately,
the result of Caviglia \cite{Caviglia2007}
 (see also Sidman \cite{Sidman2002})
 shows that 
$\reg M/I^tM = \reg (M \otimes S/I^t) \leq \reg M + \reg S/I^t = \reg M + (t-1)d+(d-1)\dim M$
where the last equality
follows from the argument above and the
 fact that $I$ is generated by a regular sequence on $S$.
\end{proof}

\begin{proof}[Proof of Proposition \ref{easy bound}]
For the first estimate we
set $d=1$ in Proposition \ref{easy bound-generalized}
and use the fact that the regularity of 
$R/I^t$ is the largest number $s$ such that
$\gm^s\not\subset I^t$. For the second estimate we
first observe that it suffices to do the case where
the number of linear forms is $\dim X$ -- that is, a fiber of the
projection is just the intersection of $X$ with a plane of complementary
dimension. Under the hypotheses given,
such a plane section of $X$ is a scheme of degree $\deg X$ and is
nondegenerate. The latter condition implies that the regularity of the
fiber is bounded above by $\deg X -\codim X+1$. Theorem \ref{main}
now gives the desired equality.
\end{proof}

\section{The Fibers of Finite Morphisms}
\label{products}

We now turn to the result that will
allow us to give the maximum regularity of the fibers of a
finite morphism in terms of global data (Corollary \ref{main1}).

\begin{theorem}
\label{main}
Let $X$ be a scheme, and let 
$\phi: X\to \PP^s$ be a finite morphism, corresponding to the 
line bundle $\cL=\phi^{*}\cO_{\PP^s}(1)$ and
the space of global sections $V=\phi^*(H^0\cO_{\PP^s}(1)) \subset H^0 \cL$.
Let $M$ be a coherent sheaf on $X$, and let
$W\subset H^0(M)$ be a space of sections. The following are
equivalent:
\begin{enumerate}
\item For every integer $t\gg 0$,
the map
$$
\Sym_t(V)\otimes W\to H^0(\cL^t\otimes M)
$$
 is surjective.
\item For every  closed point $p\in \PP^s$,
the restriction map 
$$
W\to  H^0(\phi^{-1}(p)\otimes M)
$$ 
is surjective.
\item The map of sheaves
$$
\mu: W\otimes \cO_{\PP^s} \to\phi_{*}M.
$$
is surjective
\end{enumerate}
\end{theorem}

\begin{proof} 
$1\Leftrightarrow 3$:\ \ By Serre's Vanishing Theorem, the surjectivity of $\mu$
is equivalent to the surjectivity, for $t\gg 0$, of the map
$$
W\otimes \Sym_t(V)
=W\otimes H^0(\cO_{\PP(V)}(t))
\to
 H^0(\phi_*(M)(t)).
$$
Thus
$$
\phi_*(M)(t) =\phi_*(M)\otimes_{\cO_{\PP^s}} \cO_{\PP^s}(t) = \phi_*(M\otimes_{\cO_X}\phi^*\cO_{\PP^s}(t))
=  \phi_*(M\otimes_{\cO_X}\cL^t).
$$
Taking global sections, this gives
$$
H^0(\phi_*(M)(t))
=
H^0( \phi_*(M\otimes_{\cO_X}\cL^t))=
H^0 (M\otimes_{\cO_X}\cL^t),
$$
proving that assertion 1 is equivalent to assertion 3.

$2\Leftrightarrow 3$:\ \  Since 
$\phi_*M$ is coherent, the surjectivity of $\mu$ is
equivalent by Nakayama's Lemma to the the surjectivity of all the restriction maps
$W= W\otimes \cO_{\{p\}}\to (\phi_{*}M)\otimes \cO_{\{p\}}$, where $p$ runs over
the closed points of $\PP^s$ (or just of the image of $X$).
 Using the finiteness of $\phi$, we can make
 the identifications
 \begin{eqnarray*}
(\phi_{*}M)\otimes \cO_{\{p\}}
&=
H^0( (\phi_{*}M)\otimes \cO_{\{p\}} )\\
&=
H^0( \phi_{*}(M\otimes \phi^*\cO_{\{p\}})) \\
&=
H^0( (M\otimes \phi^*\cO_{\{p\}}))\\
&=
H^0(M\otimes \cO_{\phi^{-1}(p)}),
\end{eqnarray*}
so assertion 2 is also equivalent to the surjectivity of $\mu$.
\end{proof}

\begin{corollary}\label{main1}
Suppose that $X\subset \PP^n$ is a projective scheme, and 
$\phi: X\to \PP^s$ is a finite morphism, corresponding to 
a linear system $V\subset H^0(\cL)$. The maximum regularity
of a fiber of $\phi$  over a closed point of $\PP^s$
is one more than the minimum integer $\epsilon$ such that
$$
H^0(\cO_{\PP^n}(\epsilon))\otimes \Sym_t(V) \to H^0(\cL^t(\epsilon))
$$
is surjective for $t\gg 0$.
\end{corollary}

\begin{proof}
The regularity of a fiber $\phi^{-1}(p)$ is the smallest integer $t$ such that
$H^i(\cI_{\phi^{-1}(p)}(t-i))=0$ for all $i>0$. For a non-empty fiber $Z=\phi^{-1}(p)$
 of dimension 0, only $i=1$
can be of significance, and the regularity of $Z$ is one more than
the minimum $\epsilon$ such that $H^1(\cI_{\phi^{-1}(p)})(\epsilon) = 0.$ Identifying
$\cO_Z$ with $\cO_Z(d)$, the long exact sequence in cohomology shows
that this is the least $\epsilon$ such that the restriction map
$$
H^0(\cO_{\PP^n}(\epsilon))\to H^0(\cO_Z(\epsilon)) \cong H^0(\cO_Z)
$$
is surjective. The Corollary thus follows from the equivalence 
$1 \Leftrightarrow 2$ of Theorem \ref{main} if we take
$M=\cO_{\PP^n}(\epsilon)$ and $W=H^0(\cO_{\PP^n}(\epsilon))$.
\end{proof}

\begin{proof} [Proof of Theorem \ref{geometric sample}]
In Corollary \ref{main1} take $\cL=\cO_X(1), \ M=\cO_X(e)$, 
and $W=H^0(M)$. The projection $\phi$ is finite
since the projection center does not meet $X$.  
\end{proof}

\vfill\eject

\bigskip

\vbox{\noindent Author Addresses:\par
\smallskip
\noindent{David Eisenbud}\par
\noindent{Department of Mathematics, University of California, Berkeley,
Berkeley CA 94720}\par
\noindent{eisenbud@math.berkeley.edu}\par
\smallskip
\noindent{Joe Harris}\par
\noindent{Department of Mathematics, Harvard University, Cambridge MA 02138}\par
\noindent{harris@math.harvard.edu}\par
}

\end{document}